\documentclass[12pt]{article}

\newcommand{\be}{\begin{eqnarray}}
\newcommand{\ee}{\end{eqnarray}}
\newcommand{\bee}{\begin{eqnarray*}}
\newcommand{\eee}{\end{eqnarray*}}
\newcommand{\tS}{\widetilde{S}}
\newcommand{\hS}{\widehat{S}}
\newcommand{\hW}{\widehat{W}}
\newcommand{\cF}{\mathcal{F}}
\newtheorem{teo}{Theorem}

\begin{document}

\title{A Delayed Black and Scholes Formula II }
\author{By Mercedes Arriojas\thanks {The research of this author is supported
 in part by the Center of Scientific and Human Development  
 in Venezuela, and by Southern Illinois University.}\,,  
Yaozhong Hu\thanks {The research of this author is supported in part by the
 National Science Foundation under Grant No. DMS 0204613 and No. EPS-9874732,
 matching support from the State of Kansas and General Research Fund of the
 University of Kansas.}\,,  
Salah-Eldin Mohammed\thanks {The research of this author is supported in part
 by NSF Grants DMS-9975462 and DMS-0203368, and by the University of Central
 Venezuela}\,,  \\
and Gyula Pap\thanks {The research of this author is supported in part by a 
 Fulbright fellowship.}}
 
\maketitle

\abstract{This article is a sequel to [A.H.M.P]. In [A.H.M.P], we develop an explicit formula for pricing European
 options when the underlying stock price follows a non-linear stochastic delay equation with fixed delays in the drift and diffusion terms. In this article, we look at models of the stock price described by stochastic functional differential equations with variable delays. We present a class of examples of stock dynamics with variable delays that permit an explicit form for 
the option pricing formula. As in [A.H.M.P], the market is complete with no arbitrage. This is achieved through the existence of an equivalent martingale measure. In subsequent work, the authors intend to test the models in [A.H.M.P] and the present article against real market data.}

\section{Stock price models with memory}

In this section we present models of stock price dynamics that are described by stochastic 
functional differential equations (sfde's). These models are {\it feasible}, in the sense that they admit unique solutions that are positive almost surely. 
  
Consider a stock whose price $S(t)$ at time $t$ is given by the stochastic functional  differential equation (sfde):
\begin{equation}
\left.
\begin{array}{lll}
dS(t)& = & f(t,S_t)\,dt + g(S(t-b))S(t)\,dW(t),\qquad t\in [0,T]
\nonumber \\
S(t)&=&\varphi (t), \qquad t\in [-L,0]
\end{array}
\right\}
\end{equation}
 The above sfde lives on a probability space $(\Omega, {\cal F}, P)$ with a filtration
 $({\cal F}_t)_{0 \leq t \leq T}$ satisfying the usual conditions.
The constants $L$, $b$ and $T$ are {\it positive}  with $L \geq b$.
The  space $C([-L,0],\textbf{R})$ of all continuous functions
 $\eta :[-L,0] \to \textbf{R}$ is a Banach space with the supremum norm
$$ \|\eta \|:= \sup_{s \in [-L,0]} |\eta (s)|.$$
The drift coefficient $f:[0,T]\times C([-L,0], \textbf{R}) \longrightarrow \textbf{R}$ is a given continuous functional, and $g:\textbf{R}\to \textbf{R}$ is continuous.
The initial process  $\varphi: \Omega \to C([-L,0],\textbf{R})$ is
 ${\cal F}_0$-measurable with respect to the Borel $\sigma$-algebra of
 $C([-L,0],\textbf{R})$. The process $W$ is a one-dimensional standard Brownian motion adapted to the filtration  $({\cal F}_t)_{0 \leq t \leq T}$ ; and
 $S_t \in C([-L,0],\textbf{R})$ stands for the segment
 $S_t(s):=S(t+s), \,\, s\in [-L,0],\,\, t \geq 0$. A general theory of existence and uniqueness of solutions to sfde's is provided in [Mo${}_1$] and [Mo${}_2$]. However, the results in [Mo${}_1$] and [Mo${}_2$] do not cover Hypotheses (E) below.

 Under the following hypotheses, we will demonstrate the feasibility of the model (1): That is,   
it has a unique pathwise solution such that $S(t)> 0$ almost surely for all $t \geq 0$ whenever $\varphi (t) > 0$ for all $t \in [-L,0]$.

\medskip

\noindent
{\bf{Hypotheses (E).}}
\medskip
\renewcommand{\labelenumi}{(\roman{enumi})}
\begin{enumerate}
\item There is a positive constant $L'$ such that
$$ |f (t,\eta)| \leq L' (1+\|\eta\|) $$
for all $(t,\eta) \in [0,T]\times C([-L,0],\textbf{R}).$
\item For each integer $n >0$, there is a positive constant $L_n$ such that
$$ |f(t,\eta^1)-f(t,\eta^2)| \leq L_n \|\eta^1-\eta^2\| $$
for all $(t,\eta^i) \in [0,T] \times C([-L,0], \textbf{R})$ with
$\|\eta^i\| \leq n$, i=1,2.
\item $f(t,\eta) >0$ for all $(t,\eta) \in [0,T]\times C([-L,0], \textbf{R}^+)$.
\item $g: \textbf{R} \to \textbf{R}$ is continuous.
\item $a$ and $b$ are positive constants.
\end{enumerate}
\medskip
\medskip



\begin{teo}

Assume  Hypotheses (E).
Then the sfde (1) has a pathwise unique solution $S$
 for a given ${\cal F}_0$-measurable initial process
 $\varphi: \Omega \to C([-L,0],\textbf{R})$.
Furthermore, if $\varphi (t) \geq 0$  for all $t \in [-L,0]$ a.s., then
 $S(t) \geq 0$ for all $t \geq 0$ a.s..
If in addition $\varphi (0) > 0$ a.s., then  $S(t) > 0$  for all
 $t \geq 0$ a.s..

\end{teo}

\noindent
\textsl{Proof.}


\noindent

First let $t \in [0,b]$ and let $\varphi (t) \geq 0$ a.s. for all
 $t \in [-L,0]$.
Then (1) becomes
 \begin{equation}
\left.
\begin{array}{lll}
 dS(t)&= &f(t,S_t)\,dt + g(\varphi(t-b))S(t)\,dW(t),\qquad t\in [0,b]
\nonumber \\
S(0)&=&\varphi (0).
\end{array}
\right\}
\end{equation}
Define the martingale
$$M(t):= \int_0^t g(\varphi(u-b))\, dW(u), \qquad t \in [0,b].$$
Then $S$ solves the stochastic functional differential equation (sfde)
\begin{equation}
\left.
\begin{array}{lll}
 dS(t)& = & f(t,S_t)\,dt + S(t)\,dM(t),\qquad t\in [0,b]
\nonumber \\
S_0&=&\varphi.
\end{array} \right\}
\end{equation}
Define the process $\psi:[-L,b] \times \Omega \to \textbf{R}$ as follows: 
$\psi | [0,b]$ is the solution of the linear sode
\begin{equation}
\left.
\begin{array}{lll}
d\psi(t)& = &\psi (t)\,dM(t),\qquad t\in [0,b]
\nonumber \\
\psi(0)&=& 1;
\end{array} \right\}
\end{equation}
and for all $t \in [-L,0)$, set $\psi (t) =1$. 

Define the random process $y$ to be the unique solution of the {\it random} fde
 \begin{equation}
\left.
\begin{array}{lll}
y'(t)& = & \psi (t)^{-1} f(t, \psi_t \cdot y_t),\qquad t\in [0,b]
\nonumber \\
y_0&=& \varphi.
\end{array} \right\}
\end{equation}
Observe that the above fde admits a unique global solution $y$ by virtue of the linear growth hypothesis (E)(i) and the Lipschitz condition (E)(ii).

Denote by $[M,M]$ the quadratic variation of $M$. Then, from (4), it follows that
$$\psi (t) =  \exp \{M(t)-\frac{1}{2} [M,M](t)\} > 0$$
for all $t\in [0,b]$.

Define the process $\tilde S$ by $\tilde S(t):=\psi(t)y(t)$ for $t \in [-L,b]$.
Then by the  product rule, it follows that
\begin{equation}
\left.
\begin{array}{lll}
d\tilde S(t)& = & f(t,\tilde S_t)\,dt + \tilde S(t)\,dM(t),\qquad t\in [0,b]
\nonumber \\
\tilde S_0&=&\varphi.
\end{array} \right\}
\end{equation}
Comparing (3) and (6), it follows by uniqueness that  $P$-a.s.,
 $S(t)=\tilde S(t)$  for all $ t \in [0,b]$. Now suppose that $\varphi (t) \geq 0$ a.s. for all
 $t \in [-L,0]$. Then using (5) 
and the monotonicity Hypothesis (E)(iii), it follows that  $y(t) \geq 0$ a.s. for 
all $t \in [0,l]$. If in addition $\varphi (0) > 0$   a.s., then it also follows from (5) that $y(t)> 0$  for all $t \in [0,b]$ a.s..
Hence $S(t)=\tilde S(t) > 0$ for all $t \in [0,b]$ a.s..
Using forward steps of length $b$, it is easy to see that $S(t) > 0$ a.s. for
 all $t \geq 0$.    $\diamond$

\medskip

\medskip
\noindent
{\bf{Remark.}}

\medskip
Another feasible model for the stock price is obtained by considering the sfde 
%
%
$$
\kern -0.2in\left.
\begin{array}{lll}
dS(t)& = & h(t,S^{t-a})S(t)\,dt + g(S(t-b))S(t)\,dW(t),\qquad t\in [0,T],
\\
S(t)&=&\varphi (t), \qquad t\in [-L,0].
  \end{array} \right \}  \eqno(2')
$$
with  $S^t(s):=S(t\wedge s), \,\, t,s \in [-L,T]$, and
 $h:[0,T] \times C([-L,T], \textbf{R}) \to \textbf{R}$ is a continuous
 functional.
Theorem 1 holds for the above model if Hypotheses (E) hold with 
E(iii) replaced by the following monotonicity condition:

\medskip
\noindent
(E)(iii)$'$
 For each $\xi \in C([-L,T], \textbf{R})$ with $\xi (t) \geq 0$ for all
 $t \in [-L,T]$, one has\\
\phantom{(E)(iii)$'$} $h(t, \xi) \geq 0$ for all $t \in [0,T]$.

 \medskip

The proof is analogous to that of Theorem 1.

\bigskip

\section{A stock price model with variable delay}

In this section, we give an alternative model for the stock price dynamics with
 variable delay. 
In this case we are able to develop a Black-Scholes formula for the option
 price (cf. [B.S], [Me${}_1$]).
 
Throughout this section, suppose $ h $ is a given fixed positive number.  
Denote $\lfloor t\rfloor:=kh$ if $ kh\le t< (k+1)h$. 
   
Consider a market consisting of a riskless asset $\xi$ with a variable
 (deterministic) continuous rate of return $\lambda$, and a stock $S$ satisfying the following equations   
\begin{equation}
\left.
\begin{array}{lll}
d\xi (t)& = & \lambda (t) \xi (t)\,dt    
\nonumber \\
dS(t) & = & f(t, S(\lfloor t\rfloor)) S(t) dt
            + g(t, S(\lfloor t\rfloor)) S(t) dW(t) 
\end{array}
\right\}\label{sdvde}
\end{equation}
 for $t \in (0,T]$, with initial conditions $\xi (0)=1$ and $S(0)>0$. 
The above model lives on a probability space $(\Omega, {\cal F}, P)$ with a filtration 
$({\cal F}_t)_{0 \leq t \leq T}$ satisfying the usual conditions, and a standard one-dimensional Brownian motion $W$  adapted to the filtration.  Suppose
$f :[0,T] \times \textbf{R} \to \textbf{R}$ is a continuous function.
Assume further that $g: [0,T] \times \textbf{R} \to \textbf{R}$ is continuous
 and $g(t,v) \neq 0$ for all $(t,v) \in [0,T] \times \textbf{R}$. 
 
Under the above conditions, this model is feasible: That is $S(t) > 0$ a.s. 
for all $t > 0$. This follows by an argument similar to the proof  of Theorem 1, Section 2 in [A.H.M.P]. Details are left to the reader.

Next, we will establish the completeness of the market $\{\xi(t), S(t): t \in [0,T] \}$
and the no-arbitrage property, following the approach in Section 3 in [A.H.M.P]. 
 
For $t \in [kh, (k+1)h]$, the solution of the second equation in (\ref{sdvde})
 is given by
\begin{equation}
S(t)=S (kh)\exp\left(
\int_{kh}^t g(s,S(kh))\,dW(s)+\int_{kh}^t f(s,S(kh))\,ds
-\frac12\int_{kh}^t g(s,S (kh))^2\,ds\right).
\label{e.4.5} 
\end{equation}
As in Section 3 in [A.H.M.P], let
 $$\tS(t):=\frac{S(t)}{\xi(t)}
          =S(t)e^{-\int_0^t\lambda(s)ds},\qquad t\in[0,T],$$
 be the discounted stock price process.
Again by It\^o's formula, we obtain
 \begin{eqnarray*}
  d\tS(t)&=&\frac{1}{\xi(t)}dS(t)
            +S(t)\bigg(-\frac{\lambda(t)}{\xi(t)}\bigg)\,dt\\
         &=&\tS(t)\Big[\big \{f(t,S(\lfloor t\rfloor))-\lambda(t)\big \}\,dt
                              +g(t,S(\lfloor t\rfloor))\,dW(t)\Big].
 \end{eqnarray*}
Let
 $$\hS(t):=\int_0^t \big \{f(u,S(\lfloor u\rfloor))-\lambda(u)\big \}\,du
           +\int_0^tg(u,S(\lfloor u\rfloor))\,dW(u),\qquad t\in[0,T].
$$
Then
 \begin{equation}\label{vtShS}
  d\tS(t)=\tS(t)\,d\hS(t), \qquad 0 < t < T,
 \end{equation}
 and
 \begin{equation}\label{vtS}
  \tS(t)=S(0)+\int_0^t\tS(u)\,d\hS(u),\qquad t\in[0,T].
 \end{equation}
Define the stochastic process
 $$
\Sigma (u)
   :=-\frac{\{f(u,S(\lfloor u\rfloor))-\lambda(u)\}}
           {g(u,S(\lfloor u\rfloor))},\qquad u\in[0,T].
$$
It is clear that $\Sigma (u)$ is $\cF_{\lfloor u \rfloor}^S$-measurable for each 
$u \in [0,T]$. Furthermore, by a backward conditioning argument using steps of length $h$, the reader may check that 
$$ E_P (\rho_T)=1$$
where
$$\varrho_T
   :=\exp\left\{-\int_0^T\frac{\{f(u,S(\lfloor u\rfloor))-\lambda(u)\}}
                              {g(u,S(\lfloor u\rfloor))}\,dW(u)
        -\frac{1}{2}
  \int_0^T\left|\frac{f(u,S(\lfloor u\rfloor))-\lambda(u)}
                                 {g(u,S(\lfloor u\rfloor))}\right|^2du\right\}.
 $$
(See the argument in Section 3 in [A.H.M.P] following the statement of Theorem 2.)  Hence the Girsanov theorem ([K.K], Theorem 5.5) applies, and it follows that the process
$$
 \hW(t)
  :=W(t)+\int_0^t\frac{\{f(u,S(\lfloor u\rfloor))-\lambda(u)\}}
                      {g(u,S(\lfloor u\rfloor))}du,
  \qquad t\in[0,T],
 $$
 is a standard Wiener process under the probability measure  $Q$ defined by
 $dQ :=\varrho_T\,dP$. Using (9) and the definitions of $\hS$ and $\hW$, it is easy to see that 
 $$ d\tS(t)= \tS(t)g(t,S(\lfloor t\rfloor))\,d\hW(t),\qquad t\in[0,T].$$
This implies that $\tS$ is a $Q$-martingale, and hence the market 
$\{\xi(t), S(t): t \in [0,T] \}$ has the no-arbitrage property ([K.K], Theorem 7.1).

We now establish the completeness of the market $\{\xi(t), S(t): t \in [0,T] \}$.  To do so, let $X$ be any contingent claim, viz. an integrable $\cF_T^S$-measurable non-negative random variable.  Define the process 
 $$
 M(t):=E_Q \biggl (\frac{X}{\xi(T)}\, \bigg |\,\cF_t^S \biggr )
        =E_Q \biggl (\frac{X}{\xi(T)}\,\bigg |\, \cF_t^{\hW} \biggr ),\qquad t\in[0,T].
 $$
 Then $M(t), t\in[0,T]$, is an $(\cF_t^{\hW})$-adapted $Q$-martingale. Hence, by the martingale
 representation theorem ([K.K], Theorem 9.4), there exists an $(\cF_t^{\hW})$-predictable process
 $h_1(t)$, $t\in[0,T]$, such that
 $$\int_0^T h_1(u)^2\,du<\infty\qquad a.s.,$$
 and
 $$
 M(t)= E_Q \biggl (\frac{X}{\xi(T)}\biggr ) +\int_0^th_1(u)\,d\hW(u),\qquad t\in[0,T].
 $$
Define  
 $$\pi_S(t):=\frac{h_1(t)}{\tS(t)g(t,S(\lfloor t\rfloor))},\qquad
   \pi_\xi(t):=M(t)-\pi_S(t)\tS(t),
   \qquad t\in[0,T].$$
Consider the strategy $\{(\pi_\xi(t),\pi_S(t)):t\in[0,T]\}$ which consists of
 holding $\pi_S(t)$ units of the stock and $\pi_\xi(t)$ units of the bond at
 time $t$.
The value of the portfolio at any time $t \in[0,T]$ is given by 
 $$V(t):=\pi_\xi(t)\xi(t)+\pi_S(t)S(t)=\xi(t)M(t).$$
Furthermore,
 $$
 dV(t)=\xi(t)dM(t)+M(t)d\xi(t)
        =\pi_\xi(t)d\xi(t)+\pi_S(t)dS(t), \qquad t \in (0,T].
 $$
Consequently, $\{(\pi_\xi(t),\pi_S(t)):t\in[0,T]\}$ is a self-financing
 strategy.
Clearly \ $V(T)=\xi(T)M(T)=X$. Thus the contingent claim $X$ is attainable.
This shows that the market $\{\xi(t), S(t): t \in [0,T] \}$
 is complete.

Moreover, in order for the augmented market
 $\{\xi(t), S(t), X: t \in [0,T] \}$
 to satisfy the no-arbitrage property, the price 
 of the claim $X$ must be
 $$V(t)=\frac{\xi(t)}{\xi(T)}E_Q(X\,|\,\cF_t^S)$$
at each $t \in [0,T]$ a.s.. See, e.g., [B.R] or Theorem 9.2 in [K.K].

\bigskip

The above discussion may be summarized in the following
 formula for the fair price $V(t)$ of an option on the stock whose evolution is
 described by the sdde (\ref{sdvde}).

\medskip

\begin{teo}

Suppose that the stock price $S$ is given by the sdde (\ref{sdvde}), where
 $S(0)>0$ and $g$ satisfies Hypothesis (B).
Let $T$ be the maturity time of an option (contingent claim) on the stock with
 payoff function $X$, i.e., $X$ is an $\cF_T^S$-measurable non-negative integrable
 random variable.
Then at any time $t \in [0,T]$, the fair price $V(t)$ of the option is given
 by the formula
\begin{eqnarray}
V(t)=E_Q(X \,|\,\cF_t^S)e^{-\int_t^T\lambda(s)\, ds},
\end{eqnarray}
where $Q$ denotes the probability measure on $(\Omega,\mathcal{F})$ defined by
 $dQ :=\varrho_T\,dP$ with
 $$
\varrho_t
   :=\exp\left\{-\int_0^t\frac{\{f(u,S(\lfloor u\rfloor))-\lambda(u)\}}
                              {g(u,S(\lfloor u\rfloor))}\,dW(u)
                -\frac{1}{2}
                 \int_0^t\left|\frac{f(u,S(\lfloor u\rfloor))-\lambda(u)}
                                    {g(u,S(\lfloor u\rfloor))}\right|^2
                  du\right\}
$$
 for $t\in[0,T]$.
The measure $Q$ is a martingale measure and the market is complete.

Moreover, there is an adapted and square integrable process $h_1(t),\ t\in[0, T]$, such that
\[
E_Q \biggl (\frac{X}{\xi(T)}\, \bigg |\,\cF_t^S \biggr )
 = E_Q \biggl (\frac{X}{\xi(T)}\biggr )  +\int_0^th_1(u)\,d\hW(u),\qquad t\in[0,T],
\]
where 
$$
 \hW(t)
  :=W(t)+\int_0^t\frac{\{f(u,S(\lfloor u\rfloor))-\lambda(u)\}}
                      {g(u,S(\lfloor u\rfloor))}\,du,
  \qquad t\in[0,T].
 $$
The hedging strategy is given by 
\begin{equation}
 \pi_S(t):=\frac{h_1(t)}{\tS(t)g(t,S(\lfloor t\rfloor))},\qquad
 \pi_\xi(t):=M(t)-\pi_S(t)\tS(t),
 \qquad t\in[0,T].
\end{equation}

\end{teo}

The following result gives a Black-Scholes-type formula for the value of a
 European option on the stock at any time prior to maturity ([B.S], [Me$_1$], [H.R]).

\medskip

\begin{teo}

Assume the conditions of Theorem 2. 
Let $V(t)$ be the fair price of a European call option written on the stock $S$
 with exercise price $K$ and maturity time $T$.
 Then for all $t \in \big[T-\lfloor T\rfloor,T\big]$, $V(t)$ is given by 
%
\begin{eqnarray}\label{vVBS}
V(t)=S(t)\Phi(\beta_+(t))-K\Phi(\beta_-(t))e^{-\int_t^T\lambda(s)ds}, 
\end{eqnarray}
%
 where
 $$\beta_\pm(t)
   :=\frac{\log\frac{S(t)}{K}
    +\int_t^T\left(\lambda(u)\pm\frac{1}{2}g(u,S(\lfloor u\rfloor))^2\right)du}
          {\sqrt{\int_t^Tg(u,S(\lfloor u\rfloor))^2du}},
$$
and $\Phi$ is the standard normal distribution function.\newline
If $T > h$ and $t<T-\lfloor T\rfloor$, then 
 \begin{eqnarray}\label{vVBSC}
V(t)  & =& e^{\int_0^t\lambda(s)ds}
    E_Q\left(H\left(\tS(T-\lfloor T\rfloor),
           -\frac{1}{2}\int_{T-\lfloor T\rfloor}^Tg(u,S(\lfloor u\rfloor))^2du,
\right.\right.
\nonumber\\ 
&&\qquad \qquad  \qquad \qquad \qquad \left. \left.  
                 \int_{T-\lfloor T\rfloor}^Tg(u,S(\lfloor u\rfloor))^2du\right)
             \,\bigg|\,\cF_t\right)
\end{eqnarray}
 where $H$ is given by 
 $$
H(x,m,\sigma^2):=x e^{m + \sigma^2/2} \Phi(\alpha_1(x,m,\sigma))
                   -K\Phi(\alpha_2(x,m,\sigma))e^{-\int_0^T\lambda(s)ds},
$$
 and
 $$\alpha_1 (x,m,\sigma):=\frac{1}{\sigma} \biggl [\log
\left(\frac{x}{K}\right) +\int_0^T\lambda(s)ds +m+ \sigma^2 
\biggr ],$$
$$\alpha_2 (x,m,\sigma):=\frac{1}{\sigma} \biggl [\log
\left(\frac{x}{K}\right) +\int_0^T\lambda(s)ds +m \biggr 
],$$
for $\sigma, x \in \textbf{R}^+, \, m \in \textbf{R}$.

The hedging strategy is given by
\[
 \pi_S(t)=\Phi(\beta_+(t)),\qquad
 \pi_\xi(t)=-K\Phi(\beta_-(t))e^{-\int_0^T\lambda(s)ds},\qquad
 t\in\big[T-\lfloor T\rfloor,T\big].
\]
\end{teo}

\bigskip
\bigskip

\noindent
{\bf{Acknowledgments.}}

\medskip

The authors are very grateful to R. Kuske and B. $\emptyset$ksendal  
for very useful suggestions and
corrections to earlier versions of the manuscript. The authors also
acknowledge helpful comments and discussions with Saul Jacka.

\bigskip

\centerline{\bf{ References}}

\medskip
\medskip
\par\noindent
[A.H.M.P] Arriojas, M., Hu, Y., Mohammed, S.-E.A. and  Pap, G.,
A delayed Black and Scholes formula I,
{\it preprint}, pp. 10.



\medskip
\par\noindent
[B.R] Baxter, M. and Rennie, A.,
{\it Financial Calculus}, Cambridge University Press (1996).

\medskip
\par\noindent
[B.S] Black, F. and  Scholes, M.,
The pricing of options and corporate liabilities,
{\it Journal of Political Economy} 81 (May-June 1973), 637-654.





\medskip
\par\noindent
[H.R]  Hobson, D.,  and Rogers, L. C. G.,
Complete markets with stochastic volatility,
{\it  Math. Finance\/} 8 (1998), 27--48.





\medskip
\par\noindent
[K.K] Kallianpur, G. and Karandikar R. J.,
{\it Introduction to Option Pricing Theory},
Birkh\"{a}user, Boston-Basel-Berlin (2000).



\medskip
\par\noindent
[Me${}_1$] Merton, R. C.,
Theory of rational option pricing,
{\it Bell Journal of Economics and Management Science} 4,
Spring (1973), 141-183.



\medskip
\par\noindent
[Mo${}_1$] Mohammed, S.-E. A.,
{\it Stochastic Functional Differential Equations.} Pitman 99 (1984).

\medskip
\par\noindent
[Mo${}_2$] Mohammed, S.-E. A.,
Stochastic differential systems with memory:
Theory, examples and applications. In ``Stochastic Analysis",
Decreusefond L. Gjerde J., \O ksendal B., Ustunel A.S.,  edit.,
{\it Progress in Probability} 42, Birkhauser (1998), 1-77.







\medskip
\bigskip
\noindent
Mercedes Arriojas\\
Department of Mathematics \\
University  Central of Venezuela \\
Caracas \\
Venezuela \\
e-mail: marrioja@euler.ciens.ucv.ve

\bigskip
\noindent
Yaozhong Hu \\ 
Department of Mathematics \\   
University of Kansas \\  
405 Snow Hall \\   
Lawrence \\
KS 66045-2142 \\ 
e-mail: hu@math.ukans.edu

\bigskip
\noindent
Salah-Eldin A. Mohammed \\
Department of Mathematics\\
Southern Illinois University at Carbondale\\
 Carbondale\\
 Illinois 62901, U.S.A.\\
e-mail: salah@sfde.math.siu.edu\\
web-site: http://sfde.math.siu.edu

\bigskip
\noindent
Gyula Pap \\
Institute of Mathematics and Informatics\\
University of Debrecen\\
Pf.~12, H-4010 Debrecen, Hungary\\
e-mail: papgy@inf.unideb.hu\\
web-site: http://www.inf.unideb.hu/valseg/dolgozok/papgy/papgy.html

\end{document}